\documentclass[12pt,reqno]{amsart}
\usepackage{amsfonts,amsmath,amssymb}
\usepackage{graphicx}
\usepackage{psfrag}
\usepackage{pstricks}
\usepackage{pst-tree}

\allowdisplaybreaks
\advance\textheight by \topskip

\numberwithin{equation}{section}

\def\q#1{[#1]_q}

\title[Continued fraction expansion]{A continued fraction expansion for a
	$\boldsymbol q$-tangent function: \\ An elementary proof}

\author[H.~Prodinger]{Helmut Prodinger}
\address{Helmut Prodinger\\
Department of Mathematics\\
University of Stellenbosch\\
7602 Stellenbosch\\
South Africa}
\email{hproding@sun.ac.za}

\date{May 11, 2008}

\begin{document}
\begin{abstract}
We prove a continued fraction expansion for a certain $q$-tangent
function that was conjectured by the present writer, then proved by Fulmek,
now in a completely elementary way.
\end{abstract}

\maketitle

\thispagestyle{myheadings}
\font\rms=cmr8 
\font\its=cmti8 
\font\bfs=cmbx8

\markright{\its S\'eminaire Lotharingien de
Combinatoire \bfs xx \rms (200?), Article~xxx\hfill}
\def\thepage{}


\section{Introduction}
In \cite{Prodinger}, the present writer defined the following $q$-trigonometric functions
\begin{align*}
\sin_q(z)&=\sum_{n\ge0}\frac{(-1)^nz^{2n+1}}{\q{2n+1}!}q^{n^2},
	\\
\cos_q(z)&=\sum_{n\ge0}\frac{(-1)^nz^{2n}}{\q{2n}!}q^{n^2}.
\end{align*}
Here, we use standard $q$-notation:
\begin{gather*}
\q{n}:=\frac{1-q^n}{1-q},\qquad\q{n}!:=\q{1}\q{2}\dots\q{n}.
\end{gather*}
These $q$-functions are variations of Jackson's \cite{Jackson}
$q$-sine and $q$-cosine functions.

For the $q$-tangent function $\tan_q(z)=\frac{\sin_q(z)}{\cos_q(z)}$, 
 the following continued fraction expansion was conjectured in \cite{Prodinger}:
\begin{equation*}
z\tan_q(z)
=
\cfrac{z^2}{\q{1}q^{0}-
	\cfrac{z^2}{\q{3}q^{-2}-
		\cfrac{z^2}{\q{5}q^{1}-
			\cfrac{z^2}{\q{7}q^{-9}-\dotsb
			}
		}
	}
}.
\end{equation*}
Here, the powers of $q$ are of the form ${(-1)^{n-1} n(n-1)/2-n+1}$.

In \cite{Fulmek}, this statement was proven using heavy machinery from $q$-analysis.

Happily, after about 8 years, I was now successful to provide a complete \emph{elementary}
proof that I will present in the next section.

\section{The proof}

We write
\begin{align*}
\frac {z\sin_q(z)}{\cos_q(z)}&=\cfrac{z^2}{N_0}=\cfrac{z^2}{C_1-\cfrac{z^2}{N_1}}=\cfrac{z^2}{C_1-\cfrac{z^2}{C_2-\cfrac{z^2}{N_2}}}=\dots,
\end{align*}
and set
\begin{equation*}
N_i=\frac{a_i}{b_i}.
\end{equation*}
This means
\begin{equation*}
N_i=C_{i+1}-\frac{z^2}{N_{i+1}}
\end{equation*}
or
\begin{equation*}
\frac{z^2}{N_{i+1}}=C_{i+1}-N_i
\end{equation*}
and
\begin{equation*}
\frac{b_{i+1}z^2}{a_{i+1}}=C_{i+1}-\frac{a_i}{b_i}=\frac{C_{i+1}b_i-a_i}{b_i}.
\end{equation*}
Therefore $a_i=b_{i-1}$ and
\begin{equation*}
b_{i+1}z^2=C_{i+1}b_i-b_{i-1}.
\end{equation*}
The initial conditions are
\begin{equation*}
b_{-1}=\cos_q(z)
\quad\text{and}\quad
b_0=\sum_{n\ge0}\frac{(-1)^nq^{n^2}z^{2n}}{\q{2n+1}!}.
\end{equation*}
The constants $C_i$ guarantee that all the $b_i$ are power series, i.e., they make the constant term
in $C_{i+1}b_i-b_{i-1}$ disappear. Our goal is to show that
$C_i=q^{(-1)^{i-1} i(i-1)/2-i+1}$ are the (unique) numbers that do this. We are proving the claim by proving the following \emph{explicit} formula for $b_i$:
\begin{equation*}
b_i=\sum_{n\ge0}\frac{(-1)^nz^{2n}}{\q{2n+2i+1}!}\bigg(\prod_{j=1}^{i}\q{2n+2j}\bigg)
q^{(n+\lfloor\frac{i+1}{2}\rfloor)^2+[i \text{ odd}]\binom{i+1}2}.
\end{equation*}
Note that the $C_i$ are uniquely determined by the imposed condition, and since the $b_i$ are
power series, we are done once we prove this formula by induction. The first two instances satisfy this, and we do the induction step now:
\begin{align*}
C_{i+1}&b_i-b_{i-1}=\q{2i+1}q^{(-1)^i\binom{i+1}{2}-i}\sum_{n\ge0}\frac{(-1)^nz^{2n}}{\q{2n+2i+1}!}\bigg(\prod_{j=1}^{i}\q{2n+2j}\bigg)
q^{(n+\lfloor\frac{i+1}{2}\rfloor)^2+[i \text{ odd}]\binom{i+1}2}\\
&\qquad-\sum_{n\ge0}\frac{(-1)^nz^{2n}}{\q{2n+2i-1}!}\bigg(\prod_{j=1}^{i-1}\q{2n+2j}\bigg)
q^{(n+\lfloor\frac{i}{2}\rfloor)^2+[i-1 \text{ odd}]\binom{i}2}\\
&=\sum_{n\ge0}\frac{(-1)^nz^{2n}}{\q{2n+2i+1}!}
\Bigg(\q{2i+1}\bigg(\prod_{j=1}^{i}\q{2n+2j}\bigg)
q^{(n+\lfloor\frac{i+1}{2}\rfloor)^2+[i \text{ odd}]\binom{i+1}2+(-1)^i\binom{i+1}{2}-i}\\
&\qquad-\q{2n+2i+1}\bigg(\prod_{j=1}^{i}\q{2n+2j}\bigg)
q^{(n+\lfloor\frac{i}{2}\rfloor)^2+[i -1\text{ odd}]\binom{i}2}\Bigg)\\
&=\frac1{1-q}\sum_{n\ge0}\frac{(-1)^nz^{2n}}{\q{2n+2i+1}!}\bigg(\prod_{j=1}^{i}\q{2n+2j}\bigg)\times\\
&\times\bigg((1-q^{2i+1})
q^{(n+\lfloor\frac{i+1}{2}\rfloor)^2+[i \text{ even}]\binom{i+1}2-i}
-(1-q^{2n+2i+1})q^{(n+\lfloor\frac{i}{2}\rfloor)^2+[i \text{ even}]\binom{i}2}\bigg).
\end{align*}
The last bracket in this expression can be simplified for $i$ even:
\begin{equation*}
-q^{(n+\frac{i}{2})^2+\binom i2+2i+1}(1-q^{2n})
\end{equation*}
and for $i$ odd:
\begin{equation*}
-q^{(n+\frac{i-1}{2})^2}(1-q^{2n}).
\end{equation*}
Putting everything together, we arrive at
\begin{align*}
C_{i+1}b_i-b_{i-1}
&=\sum_{n\ge0}\frac{(-1)^{n-1}z^{2n}}{\q{2n+2i+1}!}\bigg(\prod_{j=0}^{i}\q{2n+2j}\bigg)q^{(n+\lfloor\frac{i}{2}\rfloor)^2+[\text{$i$ even}]\binom {i+2}2}.
\end{align*}
Notice that the constant term vanishes, whence
\begin{align*}
b_{i+1}&=\sum_{n\ge1}\frac{(-1)^{n-1}z^{2n-2}}{\q{2n+2i+1}!}\bigg(\prod_{j=0}^{i}\q{2n+2j}\bigg)q^{(n+\lfloor\frac{i}{2}\rfloor)^2+[\text{$i$ even}]\binom {i+2}2}\\
&=\sum_{n\ge0}\frac{(-1)^{n}z^{2n}}{\q{2n+2(i+1)+1}!}\bigg(\prod_{j=1}^{i+1}\q{2n+2j}\bigg)q^{(n+\lfloor\frac{i+2}{2}\rfloor)^2+[\text{$i+1$ odd}]\binom {i+2}2},
\end{align*}
which is the announced formula.

\bibliographystyle{amsplain} 



\providecommand{\bysame}{\leavevmode\hbox to3em{\hrulefill}\thinspace}
\providecommand{\MR}{\relax\ifhmode\unskip\space\fi MR }
\providecommand{\MRhref}[2]{%
  \href{http://www.ams.org/mathscinet-getitem?mr=#1}{#2}
}
\providecommand{\href}[2]{#2}

\end{document}